\documentclass[titlepage,twoside,12pt]{article}
\usepackage{amssymb}
\usepackage{amsfonts}
\begin{document}

{\bf \Large Improvement on a Central Theory of PDEs} \\ \\

{\bf Elem\'{e}r E Rosinger} \\ \\
Department of Mathematics \\
and Applied Mathematics \\
University of Pretoria \\
Pretoria \\
0002 South Africa \\
eerosinger@hotmail.com \\ \\

{\bf Abstract} \\

One of the two basic theorems in [5] on the existence of solutions of PDEs is improved with
the use of a group analysis type argument. \\ \\

{\bf 1. Preliminaries} \\

Recently, two rather different and general, that is, {\it type independent} solution methods
have been developed for very large classes of linear and nonlinear systems of PDEs with
possibly associated initial and/or boundary conditions. One of them, [5], is using standard
Functional Analytic methods, while the other, [6,1,2,7-11], is based on a new idea in the
realms of PDEs, namely, the {\it Dedekind order completion} of spaces of smooth functions. \\

Contrary to widespread perceptions, it thus proves to be possible to implement no less than
two powerful solution methods for a very large variety of linear and nonlinear PDEs. These two
methods are {\it type independent} in the sense that they are {\it no longer} dependent on
specifics of one or another of the countless particular types of PDEs. \\
In fact, the essence of both methods is that, each in its own way is able to solve far more
general equations than PDEs. And it is precisely in this lifting to a higher level of
generality, one beyond PDEs, that the two methods attain their respective type independent
power. \\

These two solution methods have rather complementary strong, respectively, weak points. The
one in [5] does in fact deliver not only the existence of solutions, but also efficient
numerical methods for approximating them. On the other hand, the method in [6,1,2,7-11] can
deal with considerably more general equations, and among them linear and nonlinear systems of
PDEs with possibly associated initial and/or boundary conditions. \\

In this paper one of the basic theorems in [5] is improved, thus allowing for its application
to larger classes of PDEs. This theorem assumes two inequalities which prove to be {\it
sufficient} for the existence of solutions. \\

Applying group theoretic ideas, in this paper the respective to inequalities as significantly
relaxed, thus they lead to {\it weaker} sufficient conditions for the existence of
solutions. \\

As for the mentioned general mentality regarding ways of solving PDEs, we cite here two rather
typical views : \\

The 2004 edition of the Springer Universitext book "Lectures on PDEs" by V I Arnold, [3],
starts on page 1 with the statement :

\begin{quote}
"In contrast to ordinary differential equations, there is {\it no unified theory} of partial
differential equations. Some equations have their own theories, while others have no theory at
all. The reason for this complexity is a more complicated geometry ..."(italics added)
\end{quote}

The 1998 edition of the book "Partial Differential Equations" by L C Evans, [4], starts his
Examples on page 3 with the statement :

\begin{quote}
"There is no general theory known concerning the solvability of all partial differential
equations. Such a theory is {\it extremely unlikely} to exist, given the rich variety of
physical, geometric, and probabilistic phenomena which can be modelled by PDEs. Instead,
research focuses on various particular partial differential equations ..." (italics added)
\end{quote}

What other comment could one possibly make on such views, except to recall the ancient Latin
one :

\begin{quote}
"Sic transit gloria mundi ..."
\end{quote}

\bigskip
{\bf 2. The Basic Theorem} \\

For convenience, we recall here the mentioned basic theorem, see [5, p. 50, Theorem 7]. The
respective setup is as follows. We have a mapping \\

(2.1)~~~ $ F : H \longrightarrow K,~~~ F \in {\cal C}^1 $ \\

where $H$ and $K$ are two Hilbert spaces. This mapping has the property that the linear or
nonlinear system of PDEs, with possibly associated initial and/or boundary conditions, which
is to be solved can be reduced to solving in $u \in H$ the equation \\

(2.2)~~~ $ F ( u ) ~=~ 0 $ \\

and obviously, such a reduction can cover a very large variety of PDEs. \\

Now in order to solve (2.2), we associate with it the mapping $\phi : H \longrightarrow
\mathbb{R}$ given by \\

(2.3)~~~ $ \phi ( x ) ~=~ |\,|\, F ( x ) \,|\,|\,^2_K / 2,~~~ x \in H $ \\

and we take the square of the norm, in order not to affect the smoothness of $\phi$. \\

The aim is then to solve (2.2) by a least square method applied to $\phi$. The respective
result, [5, p. 50, Theorem 7], is as follows : \\

{\bf Theorem} \\

Suppose $\nabla \phi$ is locally Lipschitz and for certain $c,~ r > 0$ and $x \in H$, we
have \\

(2.4)~~~ $ |\,|\, ( \nabla \phi ) ( v ) \,|\,|\,_H ~\geq~
                       c \,|\,|\, F ( v ) \,|\,|\,_K,~~~ v \in B_r ( x ) $ \\

and \\

(2.5)~~~ $ |\,|\, F ( x ) \,|\,| ~\leq~ r\, c $ \\

where as usual, $B_r ( x )$ denotes the closed ball of radius $r$ centered at $x$. \\

Then there is $u \in B_r ( x )$, such that \\

(2.6)~~~ $ F ( u ) ~=~ 0 $ \\

{\bf Remark} \\

There is an obvious {\it conflict} between the two necessary conditions (2.4) and (2.5).
Indeed, the interest in (2.4) is in a {\it small} constant $c > 0$, while that is clearly not
convenient in (2.5). \\

Therefore, there is an interest in properly dealing with this conflict, and in this paper a
{\it group analysis} type argument is employed in this regard. \\ \\

{\bf 3. Group Analysis Type Argument} \\

Obviously, the equation (2.2) which is of our concern can be written in a large variety of
{\it equivalent} forms by subjecting it to suitable {\it transformations} of the {\it
dependent} and/or {\it independent} variables involved. And under such transformations the
sufficient conditions (2.4), (2.5) securing the existence of a solution for (2.2) may take
different forms, and specifically, the constants $c$ and $r$ involved can obtain various
values. Consequently, one or the other of these two conditions may become weaker, or possibly
stronger. And if both of them happen to become {\it weaker}, that naturally leads to a
convenient situation. \\

Our aim, therefore, is to identify transformations of the equation (2.2) which give weaker
forms of conditions (2.4), (2.5). \\

{\bf Dependent Variable Transforms.} Let us consider the case of transformations of the {\it
dependent} variable, namely, $F$. This means considering all the mappings \\

(3.1)~~~ $ A : K \longrightarrow K,~~~ A \in {\cal C}^1,~~~ A ( 0 ) ~=~ 0 $ \\

By composition with $F$ in (2.1) they give ${\cal C}^1$ mappings \\

(3.2)~~~ $ {}_AF ~=~ A \circ F : H \longrightarrow K $ \\

with the property \\

(3.4)~~~ $ u \in H,~~~ F ( u ) ~=~ 0 ~~~~\Longrightarrow~~~~ {}_AF ( u ) ~=~ 0 $ \\

Now, for the given $F$ in (2.1), let us suppose that it is of the form \\

(3.5)~~~ $ F ~=~ {}_AG ~=~ A \circ G,~~~
               \mbox{with}~~  G : H \longrightarrow K,~~~ G \in {\cal C}^1 $ \\

and the issue is to see what will the conditions (2.4), (2.5) become in terms of $G$. \\

Clearly, in order to find $G$ in (3.5) for a given $F$, it is sufficient to assume that $A$ in
(3.1) has an inverse which satisfies \\

(3.6)~~~ $ A^{-1} : K \longrightarrow K,~~~ A^{-1} \in {\cal C}^1 $ \\

In this way, the {\it group} of transformations we are interested in is given by \\

(3.7)~~~ $ {\cal A}_K ~=~ \{~ A : K \longrightarrow K ~~|~~
                      A ~~~\mbox{satisfies}~~ (3.1), (3.6) ~\} $ \\

And then, for $A \in {\cal A}_K$ and $F$ in (3.5) we have \\

(3.8)~~~ $ u \in H,~~~ F ( u ) ~=~ 0 ~~~~\Longrightarrow~~~~ G ( u ) ~=~ 0 $ \\

{\bf Independent Variable Transforms.} Alternatively, we can consider all the {\it surjective}
mappings \\

(3.9)~~~ $ B : H \longrightarrow H,~~~ B \in {\cal C}^1 $ \\

By composition with $F$ in (2.1) they give ${\cal C}^1$ mappings \\

(3.10)~~~ $ F_B ~=~ F \circ B : H \longrightarrow K $ \\

with the property \\

(3.11)~~~ $ u,~ v \in H,~~ F ( u ) ~=~ 0,~~ B ( v ) ~=~ u ~~~~\Longrightarrow~~~~
                                                                  F_B ( v ) ~=~ 0 $ \\

And now, alternatively, for the given $F$ in (2.1), let us suppose that it is of the form \\

(3.12)~~~ $ F ~=~ {}_BG ~=~ G \circ B,~~~
               \mbox{with}~~  B : H \longrightarrow H,~~~ B \in {\cal C}^1 $ \\

and then the issue is to see what will the conditions (2.4), (2.5) become in terms of $G$. \\

Here, in order to find $G$ in (3.12) for a given $F$, it is sufficient to assume that $B$ in
(3.9) has an inverse which satisfies \\

(3.13)~~~ $ B^{-1} : H \longrightarrow H,~~~ B^{-1} \in {\cal C}^1 $ \\

Therefore, this time we are interested in the {\it group} of transformations \\

(3.14)~~~ $ {\cal B}_H ~=~ \{~ B : H \longrightarrow H ~~|~~
                              B ~~~\mbox{satisfies}~~ (3.9), (3.13) ~\} $ \\

And thus for $B \in {\cal B}_H$ and $F$ in (3.12) we have \\

(3.15)~~~ $ u \in H,~~~ F ( u ) ~=~ 0 ~~~~\Longrightarrow~~~~ G ( B ( u ) ) ~=~ 0 $ \\ \\

{\bf 4. A Particular Case} \\

For a clarification of what is involved, let us consider the simplest case when $H = K =
\mathbb{R}$. Then for $x \in \mathbb{R}$, we have, see (2.3) \\

$~~~~~~ \phi ( x ) ~=~ ( F ( x ) )^2 / 2,~~~
                        \nabla \phi ( x ) ~=~ F\,' ( x ) $ \\

hence (2.4), (2.5) become \\

(4.1)~~~ $ |\, F\,' ( v ) \,| ~\geq~
                               c \,|\, F ( v ) \,|,~~~ v \in B_r ( x ) $ \\

and \\

(4.2)~~~ $ |\, F ( x ) \,| ~\leq~ r\, c $ \\

{\bf Dependent Variable Transform.} Clearly $A \in {\cal A}_K$, if and only if $A \in {\cal
C}^1 ( \mathbb{R}, \mathbb{R} )$, $A ( 0 ) = 0$ and $A$ is strictly monotonous. \\

Now let us assume (3.5) for some $A \in {\cal A}_K$, then (4.1), (4.2) become \\

(4.3)~~~ $ |\, A\,' ( G ( v ) )\, G\,' ( v ) \,| ~\geq~
                              c \,|\, A ( G ( v ) ) \,|,~~~ v \in B_r ( x ) $ \\

(4.4)~~~ $ |\, A ( G ( x ) ) \,| ~\leq~ r\, c $ \\

On the other hand, the version of (2.4), (2.5) for $G$ would be \\

(4.5)~~~ $ |\, G\,' ( v ) \,| ~\geq~
                               c_G \,|\, G ( v ) \,|,~~~ v \in B_{r_G} ( x ) $ \\

(4.6)~~~ $ |\, G ( x ) \,| ~\leq~ r_G\, c_G $ \\

for suitable $c_G,~ r_G > 0$. \\

And then the problem is to find $A$ and $G$ in (3.5), such that \\

(4.7)~~~~~~ $ (~~ (4.3), (4.4) ~~) ~~~~\Longrightarrow~~~ (~~ (4.5), (4.6) ~~) $ \\

We note that (4.3) is equivalent with \\

$~~~~~~ |\, G\,' ( v ) \,| ~\geq~
    c \,|\, A ( G ( v ) ) \,/\, A\,' ( G ( v ) ) \,|,~~~ v \in B_r ( x ) $ \\

hence if \\

(4.8)~~~ $ r_G ~\leq~ r $ \\

and \\

(4.9)~~~ $ |\, A ( G ( v ) ) \,/\, ( G ( v )\, A\,' ( G ( v ) ) ) \,|
                                   ~\geq~ c_G / c,~~~ v \in B_r ( x ) $ \\

then (4.3) implies (4.5). \\

{\bf Independent Variable Transform.} Obviously $B \in {\cal B}_H$, if and only if $B \in
{\cal C}^1 ( \mathbb{R}, \mathbb{R} )$ and $B$ is strictly monotonous. \\

Let now assume (3.12) for some $B \in {\cal B}_H$, then (4.1), (4.2) take the form \\

(4.10)~~~ $ |\, G\,' ( B ( v ) )\, B\,' ( v ) \,| ~\geq~
                     c\, |\, G ( B ( v ) ) \,|,~~~ v \in B_r ( x ) $ \\

(4.11)~~~ $ |\, G ( B ( x ) ) \,| ~\leq~ r\, c $ \\

On the other hand, the version of (2.4), (2.5) for $G$ would be \\

(4.12)~~~ $ |\, G\,' ( v ) \,| ~\geq~
                               c_G \,|\, G ( v ) \,|,~~~ v \in B_{r_G} ( x ) $ \\

(4.13)~~~ $ |\, G ( x ) \,| ~\leq~ r_G\, c_G $ \\

for suitable $c_G,~ r_G > 0$. \\

Thus the problem is to find $B$ and $G$ in (3.12), such that \\

(4.14)~~~~~~ $ (~~ (4.10), (4.11) ~~) ~~~~\Longrightarrow~~~ (~~ (4.12), (4.13) ~~) $ \\

However \\

$~~~~~~ B\,' ( v ) ~=~ 1 / ( ( B^{-1} )' \,( B ( v ) ) ),~~~ v \in H $ \\

hence (4.10) is equivalent with \\

$~~~~~~ |\, G\,' ( B ( v ) ) \,| ~\geq~
           ( c\, ( B^{-1} )' \,( B ( v ) ) ) |\, G ( B ( v ) ) \,|,~~~ v \in B_r ( x ) $ \\

thus if \\

(4.15)~~~ $ r_G ~\leq~ r $ \\

and \\

(4.16)~~~ $ ( B^{-1} )' \,( B ( v ) ) ~\geq~ c_G / c,~~~ v \in B_r ( x ) $ \\

then (4.10) implies (4.12). \\

Here we note that it is far more easy to obtain $B$ from condition (4.16), than it is to
obtain $A$ from condition (4.9). \\ \\

{\bf 5. A Simple Example} \\

Let us illustrate the above in the case of the equation, see (2.1), (2.2) \\

(5.1)~~~ $ F ( u ) ~=~ \lambda u^2 - 1 ~=~ 0,~~~ u \in \mathbb{R} $ \\

where $\lambda \in \mathbb{R}$ is given, and thus $F : \mathbb{R} \longrightarrow \mathbb{R}$.
Then (2.3) becomes \\

(5.2)~~~ $ \phi ( x ) ~=~ ( \lambda x^2 - 1 )^2 / 2,~~~ x \in \mathbb{R} $ \\

hence (4.1), (4.2) take the form \\

(5.3)~~~ $ |\, 2 \lambda v ( \lambda v^2 - 1 ) \,| ~\geq~
                           c\, |\, \lambda v^2 - 1 \,|,~~~ v \in B_r ( x ) $ \\

(5.4)~~~ $ |\, \lambda x^2 - 1 \,| ~\leq~ r\, c $ \\

We note that, if $\lambda v^2 - 1 \neq 0$, then (5.3) is equivalent with \\

$~~~~~~ 2 \,|\, \lambda v \,| ~\geq~ c $ \\

thus for given $x \in \mathbb{R}$ and $r > 0$, the {\it largest} possible $c$ in (5.3) is \\

(5.5)~~~ $ c_{\lambda,\,x,\,r} ~=~ 2 \,|\, \lambda \,|~
                                   \begin{array}{|l}
                                         ~~ 0 ~~~\mbox{if}~~ x - r \leq 0 \leq x + r \\ \\
                                         ~~ x - r ~~~\mbox{if}~~ 0 \leq x - r \\ \\
                                         ~~ - x - r ~~~\mbox{if}~~ x + r \leq 0
                                    \end{array} $ \\

It follows that in such a case (5.4) becomes \\

(5.6)~~~ $ |\, \lambda x^2 - 1 \,| ~\leq~ r\, c_{\lambda,\,x,\,r} $ \\

which is a rather difficult {\it implicit} relation in all it variables $\lambda, x$ and
$r$. \\

Let us now take $B \in {\cal B}_H$ as \\

(5.7)~~~ $ B ( v ) ~=~ \mu v,~~~ v \in \mathbb{R} $ \\

where $\mu \in \mathbb{R},~ \mu \neq 0$ is fixed, and assume (3.12). Then \\

(5.8)~~~ $ G ( v ) ~=~ ( \lambda / \mu^2 ) v^2 - 1,~~~ v \in \mathbb{R} $ \\

hence the {\it largest} possible $c$ for the version of (5.3) corresponding to $G$ is, see
(5.5) \\

(5.9)~~~ $ c_{\mu,\,\lambda,\,x,\,r} ~=~ 2 \,|\, \lambda \ / \mu^2 \,|~
                                   \begin{array}{|l}
                                         ~~ 0 ~~~\mbox{if}~~ x - r \leq 0 \leq x + r \\ \\
                                         ~~ x - r ~~~\mbox{if}~~ 0 \leq x - r \\ \\
                                         ~~ - x - r ~~~\mbox{if}~~ x + r \leq 0
                                    \end{array} $ \\

while (5.6) turns to \\

(5.10)~~~ $ |\, ( \lambda / \mu^2 ) x^2 - 1 \,| ~\leq~ r\, c_{\mu,\,\lambda,\,x,\,r} $ \\

or in view of (5.5) and (5.9), to the equivalent relation \\

(5.11)~~~ $ |\, \lambda x^2 - \mu^2 \,| ~\leq~ r\, c_{\lambda,\,x,\,r} $ \\

And now we can recapitulate. \\

With $F$ in (5.1), we had the sufficient conditions (2.4), (2.5) expressed in (5.3), (5.4),
and they came down to (5.6), with the largest possible $c_{\lambda,\,x,\,r}$ given in
(5.5). \\

Here we have to note the following {\it conflict}, see Remark in section 2 :

\begin{itemize}

\item in satisfying (5.3), there is an interest in a {\it small} $c$,

\end{itemize}

while on the other hand,

\begin{itemize}

\item in satisfying (5.4), a contrary interest appears.

\end{itemize}

Applying now the {\it group transformation} (5.7), instead of $F$ in (5.1), we obtain $G$ in
(5.8). And then the corresponding {\it transformed} version of (5.3), (5.4) comes down to, see
(5.11) \\

(5.12)~~~ $ |\, \lambda x^2 - \mu^2 \,| ~\leq~ r\, c_{\lambda,\,x,\,r} $ \\

In this way, the mentioned {\it conflict} in optimally handling the two basic necessary
conditions (5.3), (5.4) can now be approached through the {\it latitude} obtained in (5.12) by
the possibility of choosing $\mu \in \mathbb{R},~ \mu \neq 0$ arbitrarily, when compared with
(5.6), where $\mu$ does not appear. \\ \\

{\bf Conclusion} \\

The above difference between the following two relations, see (5.6), (5.12) \\

$~~~~~~ |\, \lambda x^2 - 1 \,| ~\leq~ r\, c_{\lambda,\,x,\,r} $ \\

$~~~~~~ |\, \lambda x^2 - \mu^2 \,| ~\leq~ r\, c_{\lambda,\,x,\,r} $ \\

where $\mu \in \mathbb{R},~ \mu \neq 0$ is the group parameter, illustrates the {\it existence}
of possibilities for a {\it group analysis} type argument even in that simple example of
equation (5.1) and of that simplest group transformation of the independent variable in
(5.7). \\

Within the general case of the Theorem in section 2, [5, p. 50, Theorem 7], and the consequent
attempt for an {\it optimal} approach to the {\it conflict} inherent between the two
sufficient conditions (2.4), (2.5) for the existence of solutions of very large classes of
linear and nonlinear systems of PDEs with possibly associated initial and/or boundary
conditions, the simple result in this regard in section 5 indicates the possibilities which
may exist in general. \\
The more detailed exploration of such possibilities of group analysis, applied both to
dependent and independent variables, in order to properly deal with the sufficient conditions
in the mentioned theorem, is to be presented in subsequent papers.

\end{document}